\documentclass[12pt, twoside, leqno]{article}
\usepackage{amsmath,amsthm,amscd, amsfonts, amssymb, graphicx, color}
\usepackage{amssymb,latexsym}
\usepackage{enumerate}
\usepackage{hyperref}
\usepackage{amscd}
\usepackage[all,cmtip]{xy}
\newtheorem{theorem}{Theorem}[section]
\newtheorem{lemma}[theorem]{Lemma}
\newtheorem{proposition}[theorem]{Proposition}

\theoremstyle{definition}
\newtheorem{definition}[theorem]{Definition}
\newtheorem{example}{Example}

\newtheorem{remark}{Remark}
\numberwithin{equation}{section}
\begin{document}
\baselineskip=17pt

\title{On two refinements of the bounded weak approximate identities}

\author{Mohammad Fozouni$^{1}$\footnote{Corresponding author} \ and Raziyeh Farrokhzad$^{2}$\\
\small{$^{1,2}$Department of Mathematics, Faculty of Basic Sciences and Engineering,}\\
 \small{Gonbad Kavous University,}\\
\small{P. O. Box 163, Gonbad-e Kavous, Golestan, Iran,}\\
\small{\textit{E-mails: \url{fozouni@gonbad.ac.ir} : \url{r\_farokhzad@gonbad.ac.ir}}}}
\date{}
\maketitle

\begin{abstract}
Let $A$ be a commutative Banach algebra with non-empty character space $\Delta(A)$. In this paper, we  change the concepts of convergence and boundedness in the classical notion of  bounded approximate identity. This work give us a new kind of approximate identity between bounded approximate identity and bounded weak approximate identity.  More precisely, a net $\{e_{\alpha}\}$ in $A$ is a \emph{c-w approximate identity} if for each $a\in A$, the Gel'fand transform of $e_{\alpha}a$ tends to the Gel'fand transform of $a$ in the compact-open topology and we say $\{e_{\alpha}\}$ is \emph{weakly bounded} if the image of $\{e_{\alpha}\}$ under the Gel'fand transform is bounded in $C_{0}(\Delta(A))$.

\vspace{20pt}
\end{abstract}
\section{Introduction}
In 1977, Jones and Lahr introduced a new notion of an approximate identity for a commutative  Banach algebra $A$ called "bounded weak approximate identity (b.w.a.i.)" and gave an example of a semisimple and commutative semigroup algebra with a b.w.a.i. which has no approximate identity, bounded or unbounded. Indeed, a bounded net $\{e_{\alpha}\}$ in $A$ is a b.w.a.i., if for each $a\in A$ and $\phi\in \Delta(A)$, $|\widehat{ae_{\alpha}}(\phi)-\widehat{a}(\phi)|\longrightarrow 0$. Since $\Delta(A)$ is a locally compact Hausdoroff space, one may ask the following question:

 What happens if the net $\{e_{\alpha}\}$ satisfies the following condition(s);

For each $a\in A$ and compact subset $K$ of $\Delta(A)$, $\sup_{\phi\in K}|\widehat{ae_{\alpha}}(\phi)-\widehat{a}(\phi)|\longrightarrow 0$, and what happens if we borrow the boundedness of $\{e_{\alpha}\}$ from $C_{0}(\Delta(A))$ in the sense that there exists an $M>0$ such that for each $\alpha$, $\|\widehat{e_{\alpha}}\|_{\Delta(A)}<M$?

In this paper, we investigate the above  conditions for commutative Banach algebras with non-empty character spaces. This study provide for  us a new kind of approximate identity between bounded approximate identity and bounded weak approximate identity. 
\section{Preliminaries}
 Throughout the paper, suppose that $A$ is a commutative Banach algebra and suppose $\Delta(A)$ is the character space of $A$, that is, the space consisting of all non-zero homomorphisms from $A$ into $\mathbb{C}$.

A bounded net $\{e_{\alpha}\}$ in $A$ is called a bounded approximate identity (b.a.i.) if for all $a\in A$, $\lim_{\alpha}ae_{\alpha}=a$.
The notion of a bounded approximate identity first arose in Harmonic analysis; see \cite[Section 2.9]{Dales} for a full discussion of approximate identity and its applications.

For $a\in A$, we define $\widehat{a}: \Delta(A)\longrightarrow \mathbb{C}$ by $\widehat{a}(\phi)=\phi(a)$ for all $\phi\in \Delta(A)$. Then $\widehat{a}$ is in $C_{0}(\Delta(A))$ and   $\widehat{a}$ is called the Gel'fand transform  of $a$. Note that $\Delta(A)$ is equipped with the Gel'fand topology which turns $\Delta(A)$ into a locally compact Hausdoroff space; see \cite[Definition 2.2.1, Theorem 2.2.3(i)]{Kaniuth2}. Since for each $\phi\in \Delta(A)$, $\|\phi\|\leq 1$; see \cite[Lemma 2.1.5]{Kaniuth2}, we have $\|\widehat{a}\|_{\Delta(A)}\leq \|a\|$, where $\| \cdot \|$ denotes the norm of $A$.

Suppose that $\phi\in \Delta(A)$. We denote by $A_{c}$ the space of all $a\in A$ such that $\mathrm{supp\ }\widehat{a}$ is compact, and by $J_{\phi}$  the space of all $a\in A_{c}$ such that $\phi\notin \mathrm{supp\ }\widehat{a}$. Also, let $M_{\phi}=\ker(\phi)=\{a\in A : \phi(a)=0\}$.

 Let $X$ be a non-empty locally compact Hausdoroff space. A subalgebra $A$ of $C_{0}(X)$ is called a function algebra if $A$ separates strongly the points of $X$, that is, for each $x,y\in X$ with $x\neq y$, there exists $f\in A$ such that $f(x)\neq f(y)$, and for each $x\in X$, there exists $f\in A$ with $f(x)\neq 0$.
   A function algebra $A$ is called a Banach function algebra if $A$ has a norm $\|\cdot\|$ such that  $(A, \|\cdot\|)$ is a Banach algebra. A Banach function algebra $A$ is called natural if $X=\Delta(A)$, that is, every character of $A$ is an evaluation functional on some $x\in X$, or $x\longrightarrow \phi_{x}$ is a homeomorphism. In the latter case $M_{\phi_{x}}$ is denoted by $M_{x}$. A uniform algebra on $X$ is a Banach function algebra $(A, \|\cdot\|)$ for which $\|\cdot\|$ is equivalent to the uniform norm $\|\cdot\|_{X}$ of $C_{0}(X)$, where $\|\cdot\|_{X}$ defined by $\|f\|_{X}=\sup\{|f(x)| : x\in X\}$   for $f\in C_{0}(X)$.

 Let $S$ be a non-empty set. By $c_{00}(S)$, we denote the space of all functions on $S$ of finite support. A Banach sequence algebra on $S$ is a Banach function algebra $A$ on $S$ such that $c_{00}(S)\subseteq A$.

A net $\{e_{\alpha}\}$ in $A$ is called a bounded weak approximate identity (b.w.a.i.) if there exists a non-negative constant $C$ such that for each $\alpha$, $\|e_{\alpha}\|<C$  and for all $a\in A$ and $\phi\in\Delta(A)$
\begin{equation*}
\lim_{\alpha}|\phi(ae_{\alpha})-\phi(a)|=0,
\end{equation*}
or equivalently, $\lim_{\alpha}\phi(e_{\alpha})=1$ for each $\phi\in \Delta(A)$.
 See \cite{Jones} and \cite{LF} for more details.

 If $\Delta(A)=\emptyset$, every bounded net in $A$ is a b.w.a.i. for $A$. So, to avoid trivialities we will always assume that $A$ is a Banach algebra with $\Delta(A)\neq \emptyset$.

 In the case that $A$ is a natural Banach function algebra, a bounded weak approximate identity $\{u_{\alpha}\}$  is called a bounded pointwise approximate identity (BPAI); see \cite[Definition 2.11]{DU}.

Let $A$ be a  Banach function algebra. We say that $A$ has   bounded relative approximate units (BRAUs) of bound $m$ if for each non-empty compact subset $K$ of $\Delta(A)$ and $\epsilon>0$, there exists $f\in A$ with $\|f\|<m$ and $|1-\varphi(f)|<\epsilon$ for all $\varphi\in K$; see \cite[Section 2.2]{DU}.

In this paper,  we introduce two notions of approximate identities of a Banach algebra $A$ depending on its character space and provide some illuminating examples to show the difference of our notions from those previously known.
\section{Definitions}
 Suppose  $\mathcal{K}(\Delta(A))$ denotes the collection of all compact subsets of $\Delta(A)$, and $\tau_{co}$ denotes the compact-open topology of $C_{0}(\Delta(A))$.

\begin{definition}\label{maindef}
A  \emph{c-w approximate identity} for $A$ (in which, by $"$c-w$"$ we mean $"$compact-weak$"$ ) is a net $\{e_{\alpha}\}$ in $A$ such that for each $a\in A$ and $K\in \mathcal{K}(\Delta(A))$
\begin{equation*}\label{formula1}
\lim_{\alpha}\|\widehat{ae_{\alpha}}-\widehat{a}\|_{K}=0.
\end{equation*}
If the net $\{e_{\alpha}\}$ is bounded, we say that it is a bounded  c-w approximate identity (b.c-w.a.i.) for Banach algebra $A$.

 The Banach algebra  $A$ has  \emph{c-w approximate units} (c-w.a.u.), if  for each $a\in A$, $K\in \mathcal{K}(\Delta(A))$ and $\epsilon>0$, there exists $e\in A$  such that $\|\widehat{ae}-\widehat{a}\|_{K}<\epsilon.$
 The  approximate units have bound $m$, if $e$ can be chosen with $\|e\|<m$. Clearly, if $A$ has a b.c-w.a.i., then it has b.c-w.a.u.
  \end{definition}
 Recall that $A[\tau]$ denotes the topological algebra $A$, where $\tau$ indicates the topology of the underlying topological vector space $A$; see \cite{Frag} for a general theory of topological algebras.
\begin{definition}\label{dfn}
A net $\{e_{\lambda}\}$ in $A$ is a \emph{weakly bounded c-w approximate identity} (w.b.c-w.a.i.) if the net $\{\widehat{e_{\lambda}}\}$ is a b.a.i. for the topological algebra $\widehat{A}[\tau_{co}]$, that is, for each $a\in A$ and  $K\in \mathcal{K}(\Delta(A))$,
$\lim_{\alpha}\|\widehat{ae_{\lambda}}-\widehat{a}\|_{K}=0$ and there exists a constant $M>0$ such that for all $\lambda$ and $K\in \mathcal{K}(\Delta(A))$,
\begin{equation}\label{bd}
P_{K}(\widehat{e_{\lambda}})=\sup\{|\phi(e_{\lambda})| : \phi\in K\}=\|\widehat{e_{\lambda}}\|_{K}<M.
\end{equation}
Since, each $\widehat{e_{\lambda}}$ is in $C_{0}(\Delta(A))$,  relation (\ref{bd}) is equivalent to $\|\widehat{e_{\lambda}}\|_{\Delta(A)}< M$.

Also,  $A$ has a \emph{weakly bounded  c-w approximate units} (w.b.c-w.a.u.) of bound $m$, if  for each $a\in A$, $K\in \mathcal{K}(\Delta(A))$ and $\epsilon>0$, there exists $e\in A$ with $\|\widehat{e}\|_{\Delta(A)}<m$  and $\|\widehat{ae}-\widehat{a}\|_{K}<\epsilon.$
\end{definition}
It is a routine calculation that each b.c-w.a.i. is a w.b.c-w.a.i.
 We will show in the next section that this two concepts are different.

 For a natural uniform algebra $A$, one can see immediately that each w.b.c-w.a.i. is a b.c-w.a.i.
 \begin{remark}
 We can give a more general version of Definition \ref{dfn} as follows:

  Let $\tau$ be a topology on $C_{0}(\Delta(A))$ generated by the saturated  family $(p_{i})$ of seminorms on $\Delta(A)$; see \cite[ Definition 1.7]{Frag} for more details of locally convex topological algebras and saturated family. A net $\{a_{\alpha}\}$ in $A$ is a w.b.c-w.a.i. for $A$ if $\{\widehat{a_{\alpha}}\}$ is a b.a.i. for the topological algebra $\widehat{A}[\tau]$, that is, for each $a\in A$, $\widehat{a_{\alpha}a}\xrightarrow{\tau} \widehat{a}$ and there exists $M>0$ such that for each $i$ and $\alpha$, $p_{i}(\widehat{a_{\alpha}})<M$.
   But in the current paper we only focus on the case that $\tau=\tau_{co}$.
 \end{remark}
\section{Examples}

For a locally compact group $G$ and $1<p<\infty$, let $A_{p}(G)$ denote the Fig$\mathrm{\grave{a}}$ Talamanca-Herz algebra which is a natural Banach function algebra on $G$; see \cite{Herz}. The group $G$ is said to be amenable if there exists an $m\in L^{\infty}(G)^{*}$ such that $m\geq 0$, $m(1)=1$ and $m(L_{x}f)=m(f)$ for each $x\in G$ and $f\in L^{\infty}(G)$,  where $L_{x}f(y)=f(x^{-1}y)$; see \cite[Definition 4.2]{Pier}.  A classical theorem due to Leptin and Herz, characterize the amenability of a group $G$ through the existence of a bounded approximate identity for the  Fig$\mathrm{\grave{a}}$ Talamanca-Herz algebra.

The following example provide for us  a Banach algebra with a w.b.c-w.a.i. such that has no  b.c-w.a.i.
\begin{example}\label{ex}
Let $1<p<\infty$ and $G$ be a non-amenable locally compact group. Using \cite[Proposition 3.11]{DU} we can see that $A_{p}(G)$ does not have any b.c-w.a.i. Now, we construct a w.b.c-w.a.i. for $A_{p}(G)$.
Put $\Lambda=\{K\subseteq G : \textrm{$K$ is compact and } |K|>0\}$. It is obvious that $\Lambda$ with inclusion is a directed set.
For each $K\in \Lambda$ define $u_{K}$ as follows,
\begin{equation*}
u_{K}:=|K|^{-1}\chi_{KK}\ast \check{\chi}_{K}.
\end{equation*}
Clearly, $\{u_{K}\}$ is a net in $A_{p}(G)$. For each $x\in G$ we have
\begin{align*}
u_{K}(x)=|K|^{-1}\int_{G}\chi_{KK}(y)\check{\chi}_{K}(y^{-1}x)dy&=|K|^{-1}\int_{KK}\chi_{K}(x^{-1}y)dy\\
&=|K|^{-1}\int_{KK}\chi_{xK}(y)dy\\
&=\frac{|KK\cap xK|}{|K|}.
\end{align*}
If $x\in K$, $KK\cap xK=xK$. Therefore, $u_{K}(x)=1$ and otherwise since $KK\cap xK\subseteq xK$, $0\leq u_{K}(x)\leq 1$.
Hence, the net $\{\widehat{u_{K}}\}$ is  bounded  in $C_{0}(G)$.

Now, let $f$ be an arbitrary element of $A_{p}(G)$ and $K^{'}$ be a compact subset of $G$. Since $G$ is a locally compact group, for each $x\in K^{'}$ there exists a compact neighborhood $V_{x}$ of $x$. On the other hand, we know that $K^{'}\subseteq \cup_{x\in K^{'}}V_{x}$ and for each $x$, $|V_{x}|>0$. But $K^{'}$ is compact, so there are points $x_{1},\ldots, x_{n}$ in $K^{'}$ such that $K^{'}\subseteq \cup_{i=1}^{n}V_{x_{i}}$. Therefore, by putting $K^{''}=\cup_{i=1}^{n}V_{x_{i}}$ we conclude that $K^{''}\in\Lambda$  and $K^{'}\subseteq K^{''}$.
Now, it is obvious that $\lim_{K\in \Lambda}\|\widehat{u_{K}f}-\widehat{f}\|_{K^{'}}=0$ and this completes the proof.
\end{example}
Let $G$ be a locally compact group, $A(G)$ be the Fourier algebra, and $L^{1}(G)$ be the group algebra endowed with the norm $\|\cdot\|_{1}$ and the convolution product. Put $\mathfrak{L}A(G)=L^{1}(G)\cap A(G)$ with the norm $|||f|||=\|f\|_{1}+\|f\|_{A(G)}$.
We know that $\mathfrak{L}A(G)$ with the pointwise multiplication is a commutative Banach algebra   called  the Lebesgue-Fourier algebra of $G$ and $\Delta(\mathfrak{L}A(G))=G$; see \cite{GHL}.
It was shown that $\mathfrak{L}A(G)$ has a b.a.i. if and only if $G$ is a compact group; see \cite[Proposition 2.6]{GHL}.

The following example provide for us  a Banach algebra with a w.b.c-w.a.i whereas has no  b.a.i.
\begin{example}
Let $G=\mathbb{R}$ be the real line additive group and $A=\mathfrak{L}A(G)$. Clearly $G$ is not compact and hence $A$ has no b.a.i. On the other hand, it is well-known that $G$ is amenable, and hence $A(G)$ has a b.a.i. $\{u_{\alpha}\}$ in $A(G)\cap C_{c}(G)$ by the Leptin-Herz theorem. Also, it is well-known that for each $u\in A(G)$, $\|u\|_{G}\leq \|u\|_{A(G)}$.
 So, for each $u\in A$ and $K\in \mathcal{K}(G)$, we have
$$\sup_{x\in K}|u(x)u_{\alpha}(x)-u(x)|\leq \sup_{x\in G}|u(x)u_{\alpha}(x)-u(x)|\leq \|uu_{\alpha}-u\|_{A(G)}\longrightarrow 0.$$
Therefore, $\{u_{\alpha}\}$ is a w.b.c-w.a.i. for $A$.
\end{example}

Recall that if $A$ is a function algebra on $K$, then $x\in K$ is a peak point  if there exists $f\in A$ such that $f(x)=1$ and $|f(y)|<1$ for each $y\in K\setminus\{x\}$. It is well-known that for the disc algebra $A(\overline{\mathbb{D}})$; see \cite[Example 2.1.13(ii)]{Dales} for more details on disc algebra,  $z\in \overline{\mathbb{D}}$ is a peak point if and only if $z\in \mathbb{T}$; $\mathbb{T}=\{z\in \mathbb{C} : |z|=1\}$ and $\mathbb{D}=\{z\in \mathbb{C} : |z|< 1\}$.

It is worth noting that there exist Banach algebras without any w.b.c-w.a.i. as the following example shows.

\begin{example}
Let $A=A(\overline{\mathbb{D}})$ be the disc algebra and for $z_{0}\in \mathbb{ D}$, let $B=M_{z_{0}}$. Clearly, $\overline{\mathbb{D}}\setminus \{z_{0}\}\subseteq \Delta(B)$. So, if $B$ has a w.b.c-w.a.i., then $B$ has a BPAI which is in contradiction with \cite[Example 4.8(i)]{DU}.
\end{example}
Note that if a Banach function algebra $A$ has a b.c-w.a.i., then it has BRAUs. The following example gives a Banach function algebra with a BPAI, while it does not have any BRAUs. So, two concepts of BPAI and b.c-w.a.i. are different.
\begin{example}
Let $\mathbb{I}=[0,1]$ and let $A=\{f\in C(\mathbb{I}) : I(f)<\infty\}$, where $$I(f)=\int_{0}^{1}\frac{|f(t)-f(0)|}{t}dt.$$
For each $f\in A$, define $\|f\|=\|f\|_{\mathbb{I}}+I(f)$. By \cite[Example 5.1]{DU}, $(A, \|\cdot\|)$ is a natural Banach function algebra. Also, $M_{0}$ does not have BRAUs, but it has a BPAI. Therefore, $M_{0}$ does not have any b.c-w.a.i.
\end{example}
The following example shows the difference between b.a.i. and b.c-w.a.i.
\begin{example}\label{Ex: 1}
Let $\alpha=(\alpha_{k})\in \mathbb{C}^{\mathbb{N}}$ and for each $n\in \mathbb{N}$,  set
\begin{equation*}
p_{n}(\alpha)=\frac{1}{n}\sum_{k=1}^{n}k|\alpha_{k+1}-\alpha_{k}|,\quad p(\alpha)=\sup_{n\in \mathbb{N}}p_{n}(\alpha).
\end{equation*}
Put $A=\{\alpha \in c_{0} : p(\alpha)<\infty\}$. By \cite[Example 4.1.46]{Dales}, $A$ is a natural Banach sequence algebra on $\mathbb{N}$ for the norm given by $\|\alpha\|=\|\alpha\|_{\mathbb{N}}+p(\alpha)$. For each $K\in \mathcal{K}(\mathbb{N})$ there exists $\alpha_{K}\in A$ such that $\alpha_{K}(k)=1$ for all $k\in K$ and $\|\alpha_{K}\|\leq 4$. So, $\{\alpha_{K} : K\in \mathcal{K}(\mathbb{N})\}$ is a b.c-w.a.i. for $A$. But by parts (iii) and (v) of \cite[Example 4.1.46]{Dales}, $A^{2}$ has  infinite codimension in $A$ where $A^{2}=\mathrm{linear\ span}\{ab : a,b\in A\}
$. Therefore, $A$ has no b.a.i.
\end{example}
\begin{remark}
By Cohen's factorization theorem, we know that if a Banach algebra $A$ has a b.a.i., then $A$ factors, that is, for each $a\in A$ there exists $b, c\in A$ such that $a=bc$. One  may ask this question: Whether the statement of the Cohen theorem is valid if we replace b.a.i. by b.c-w.a.i. (w.b.c-w.a.i.)?  Example \ref{Ex: 1} gives a negative answer to this question. Indeed, if $A$ factors, then $A^{2}=A$ which is a contradiction.
\end{remark}
Let $G$ be a locally compact abelian group and let $A=S(G)$ be a Segal algebra in $L^{1}(G)$; see \cite[Definition 4.5.26]{Dales} or \cite{Reiter}.
By using the Cohen factorization theorem, if $A$ has a b.a.i., then $A=L^{1}(G)$. But with the aid of  \v{S}ilov's idempotent theorem we show in the sequel  that there exist a Segal algebra $S(G)$  with a w.b.c-w.a.i. such that does not satisfy $S(G)=L^{1}(G)$.

First we recall the \v{S}ilov idempotent theorem as follows; see \cite[Theorem 3.5.1]{Kaniuth2} or \cite[Theorem 2.4.33]{Dales}.
\begin{theorem}
Let $A$ be a commutative Banach algebra and $C$ be a compact and  open subset of $\Delta(A)$. Then there exists an idempotent $a\in A$ such that $\widehat{a}$ is equal to the characteristic function of $C$.
\end{theorem}
\begin{proposition}\label{cor}
Every commutative Banach algebra $A$ with discrete character space has a w.b.c-w.a.i. But the converse is not valid in general.
\end{proposition}
\begin{proof}
Let $\mathcal{F}(\Delta(A))$ be the collection of all finite subsets of $\Delta(A)$ and let $K$ be an element of $\mathcal{F}(\Delta(A))$.  So, by using \v{S}ilov's idempotent theorem, we can take an element $e_{K}$ in $A$ such that $\widehat{e_{K}}=\chi_{K}$. So, for each $K\in \mathcal{F}(\Delta(A))$, $\|\widehat{e_{K}}\|_{\Delta(A)}=1$.

Now, it is clear that $\{e_{K} : K\in \mathcal{F}(\Delta(A))\}$ is a w.b.c-w.a.i. for $A$, where $\{K : K\in \mathcal{F}(\Delta(A))\}$ is ordered with inclusion. Because for each $a\in A$ and $K^{'}\in \mathcal{F}(\Delta(A))$,
$\lim_{K}\|\widehat{e_{K}a}-\widehat{a}\|_{K^{'}}=0.$

To see that the converse is not valid consider $A_{p}(G)$, where $G=SL(2,\mathbb{R})$ is the  multiplicative group of all $2 \times 2$ real matrices with determinant 1. We know that $G$ is non-amenable; see \cite[Exercise 1.2.6 (viii)]{Runde}, and by \cite[Proposition 1.4, pp. 207]{Sugi}, $G$ is a connected group and hence it is not discrete. By Example \ref{ex}, $A_{p}(G)$ has a w.b.c-w.a.i., but $\Delta(A_{p}(G))=G$ is not discrete.
\end{proof}
\begin{example}
Let $G$ be a locally compact abelian group with dual group $\widehat{G}$. For each $1<p\leq \infty$,  put $S_{p}(G)=L^{1}(G)\cap L^{p}(G)$ and define a norm as
$$\|f\|_{S_{p}(G)}=\max\{\|f\|_{1}, \|f\|_{p}\}\quad(f\in S_{p}(G)).$$
Then $S_{p}(G)$ is a Segal algebra. If $G$ is  a compact and infinite group, then by \cite[Remark 2]{IT}, $S_{p}(G)$ has no b.w.a.i.
Let $A=S_{p}(G)$. We know that $\Delta(A)$ is homeomorphic to $\widehat{G}$; see \cite{Reiter}. Since $G$ is a compact group, it is well-known that  $\widehat{G}$ is discrete. So, as an application of Proposition \ref{cor}, we see that $A$ has a w.b.c-w.a.i.
\end{example}


\section{Hereditary Properties}
In this section we will show that for some certain closed ideals of a Banach algebra $A$ with some conditions, $I$ has a b.c-w.a.i. (w.b.c-w.a.i.)  if and only if $A$ has a b.c-w.a.i. (w.b.c-w.a.i.). First, we give the following result which shows the relation between b.c-w.a.u. (w.b.c-w.a.u.)  and b.c-w.a.i. (w.b.c-w.a.i.),  and it is a key tool in the sequel.
\begin{proposition}\label{prop}
Let $A$ be a commutative Banach algebra. Then $A$ has a b.c-w.a.u. (w.b.c-w.a.u.) if and only if $A$ has a b.c-w.a.i. (w.b.c-w.a.i.).
\end{proposition}
\begin{proof}
Let $A$ has a b.c-w.a.u. of bound $M>0$. Suppose that $\mathcal{F}$ is a finite subset of $A$, $K\in \mathcal{K}(\Delta(A))$ and $\epsilon>0$. Then like the proof of \cite[Proposition 1.1.11]{Kaniuth2}, there exists $e_{(\mathcal{F},K,\epsilon)}\in A$ such that $\|e_{(\mathcal{F},K,\epsilon)}\|< M$ and
$$\|\widehat{ae_{(\mathcal{F},K,\epsilon)}}-\widehat{a}\|_{K}<\epsilon\quad (a\in \mathcal{F}).$$
Now, consider the following net, $$\mathfrak{U}=\{e_{(\mathcal{F},K,1/n)} : K\in \mathcal{K}(\Delta(A)), \mathcal{F}\subseteq A \text{ is finite and } n\in \mathbb{N}\}.$$ It is a routin calculation that $\mathfrak{U}$ is a b.c-w.a.i. for $A$.
\end{proof}
Let $I$ be a closed ideal of $A$. Suppose that $I$ and $A/I$; the quotient Banach algebra,
 respectively have  b.a.i. of bound $m$ and   $n$. Then $A$ has a b.a.i. of bound $m+n+mn$; see \cite[Lemma 1.4.8 (ii)]{Kaniuth2}. In the setting of b.c-w.a.i. we have the following version of the mentioned assertion.
\begin{lemma}\label{Lem}
Suppose that $I$ has a w.b.c-w.a.i. (b.c-w.a.i.) of bound $m$ and $A/I$ has a b.a.i. of bound $m$. Then $A$ has a w.b.c-w.a.i. (b.c-w.a.i.) of bound $m+n+mn$.
\end{lemma}
\begin{proof}
With a slight modification in the proof of \cite[Lemma 1.4.8 (ii)]{Kaniuth2}, we can see the proof. Therefore, we omit the details.
\end{proof}
There exists a Banach function algebra $A$ such that $A$ has a b.c-w.a.i., but one of its closed ideals has no b.c-w.a.i.; see \cite[Example 5.6]{LF}. Indeed, $A=C^{1}[0,1]$; the algebra of all functions with continuous derivation, has a b.c-w.a.i., but for each $t_{0}\in [0,1]$, $M_{t_{0}}$ has no BPAI and hence no b.c-w.a.i.  So,  the converse of  Lemma \ref{Lem} is not valid in general. Although  we have the following result.
\begin{theorem}
Let $A$ be a Banach algebra, $\phi_{0}\in \Delta(A)$ and $\overline{J_{\phi_{0}}}=M_{\phi_{0}}$. Suppose that there exists $n>0$ such that for each neighborhood $U$ of $\phi_{0}$, there exists $a\in A$ with $\phi_{0}(a)=1$, $\|a\|\leq n$ and $\mathrm{supp\ } \widehat{a}\subseteq U$.  Then $M_{\phi_{0}}$ has a w.b.c-w.a.i. (b.c-w.a.i.) if and only if $A$ has a w.b.c-w.a.i. (b.c-w.a.i).
\end{theorem}
\begin{proof} We only give the proof for the case $"$w.b.c-w.a.i.$"$, since  the proof for $"$b.c-w.a.i.$"$ is similar.

Suppose that $M_{\phi_{0}}$ has a w.b.cw-a.i. Clearly,  $M_{\phi_{0}}$ has codimension 1, that is, $A/M_{\phi_{0}}$ is generated by one vector. Therefore, $A/M_{\phi_{0}}$ has a b.a.i. and hence by Lemma \ref{Lem}, $A$ has a w.b.c-w.a.i.

Conversely, suppose that $A$ has a w.b.c-w.a.i. Let $a\in M_{\phi_{0}}$ and $\varepsilon>0$. Since $\overline{J_{\phi_{0}}}=M_{\phi_{0}}$, there exists $a_{1}\in J_{\phi_{0}}$ such that
\begin{equation*}
\|a-a_{1}\|<\varepsilon.
\end{equation*}
Therefore, for each $W\in \mathcal{K}(\Delta(A))$, there exists $m>0$ and $b\in A$ such that
\begin{equation*}
\|\widehat{a_{1}}-\widehat{a_{1}b}\|_{W}<\varepsilon,\quad \|\widehat{b}\|_{\Delta(A)}\leq m.
\end{equation*}
There exists a neighborhood $U$ of $\phi_{0}$ in $\Delta(A)$ such that $U\cap \mathrm{supp\ }\widehat{a_{1}}=\emptyset$, because $a_{1}\in J_{\phi_{0}}$ and $\Delta(A)$ is Hausdoroff. Now, by the hypothesis, there exists  $c\in A$ such that $\|c\|\leq n$, $\phi_{0}(c)=1$  and $\mathrm{supp\ }\widehat{c}\subseteq U$. Therefore, $\widehat{a_{1}}\widehat{c}=0$, $b-bc$ is in $M_{\phi_{0}}$ and $\|\widehat{b-bc}\|_{\Delta(A)}\leq m(n+1)$. So, for $W\in \mathcal{K}(\Delta(A))$ we have
\begin{align*}
\|\widehat{a}-\widehat{a}(\widehat{b-bc})\|_{W}&\leq \|\widehat{a}-\widehat{a_{1}}\|_{W}+\|\widehat{a_{1}}-\widehat{a_{1}}\widehat{b}\|_{W}+\|\widehat{a_{1}}\widehat{b}-\widehat{a}(\widehat{b-bc})\|_{W}\\
&\leq \|a-a_{1}\|+\varepsilon+\|\widehat{a_{1}}\widehat{b}-\widehat{a}(\widehat{b-bc})+\widehat{b}\widehat{a_{1}}\widehat{c}\|_{W}\\
&\leq 2\varepsilon+\varepsilon m(1+n).
\end{align*}
Hence, by Proposition \ref{prop}, $M_{\phi_{0}}$ has  a w.b.c-w.a.i., which completes the proof.
\end{proof}

\bibliographystyle{amsplain}
\bibliography{References}
\end{document}